\documentclass[11pt,reqno]{amsart}
\usepackage[T1]{fontenc}
\usepackage[utf8]{inputenc}

\usepackage{amsmath,amssymb,mathtools}
\usepackage{bm}
\usepackage{mathrsfs}

\usepackage[a4paper,margin=3cm]{geometry}
\usepackage{enumitem}

\usepackage{xcolor}
\usepackage[hidelinks]{hyperref}
\usepackage[nameinlink,noabbrev]{cleveref}

\newtheorem{theorem}{Theorem}[section]
\newtheorem{proposition}[theorem]{Proposition}

\newtheorem*{assumptionA}{Assumption A}

\theoremstyle{definition}
\newtheorem{definition}[theorem]{Definition}

\newtheorem{example}[theorem]{Example}
\usepackage{needspace}
\usepackage{etoolbox}

\preto{\section}{\Needspace{6\baselineskip}}
\preto{\section}{\Needspace{12\baselineskip}}
\preto{\subsection}{\Needspace{7\baselineskip}}
\theoremstyle{remark}
\newtheorem{remark}[theorem]{Remark}

\newcommand{\Pp}{\mathbb{P}}

\newcommand{\T}{\mathcal{T}}
\newcommand{\Sspace}{\mathcal{S}^{2}}

\newcommand{\Ef}[2]{\mathcal{E}^{f}_{#1,#2}}
\newcommand{\Snell}{\mathcal{R}^{f}}
\newcommand{\esssup}{\operatorname*{ess\,sup}}
\newcommand{\essinf}{\operatorname*{ess\,inf}}
\newcommand{\ind}{\mathbf{1}}

\title[Nonlinear Optimal Stopping and Domineering Claims]
{Nonlinear Optimal Stopping and Domineering Claims}

\author{
Siham BOUHADOU\textsuperscript{$\ast$}
\quad
Oussama LAMHADOUL\textsuperscript{$\ast$}
\quad
Youssef OUKNINE\textsuperscript{$\ast,\dagger$}
}

\makeatletter
\let\originalsetauthors\@setauthors

\renewcommand{\@setauthors}{%
    \originalsetauthors
    \begingroup
    \normalfont\footnotesize
    \raggedright

    \medskip

    \noindent\textsuperscript{$\ast$}
    Department of Mathematics,
    Faculty of Sciences Semlalia,
    Cadi Ayyad University,
    B.P. 2390, Marrakech, Morocco.\par

    \noindent
    \textit{E-mail:}
    \texttt{s.bouhadou@uca.ac.ma},
    \texttt{o.lamhadoul.ced@uca.ac.ma},
    \texttt{ouknine@uca.ac.ma}\par

    \medskip

    \noindent\textsuperscript{$\dagger$}
    Mohammed VI Polytechnic University,
    Africa Business School,
    Avenue Mohammed Ben Abdellah Regragui,
    Madinat Al Irfane,
    BP 6380, Rabat, Morocco.\par

    \noindent
    \textit{E-mail:}
    \texttt{youssef.ouknine@um6p.ma}\par

    \endgroup
}
\makeatother


\date{}

\begin{document}

\begin{abstract}
We study nonlinear optimal stopping and domineering
claims in a Brownian--Poisson filtration, under a
nonlinear evaluation induced by a Lipschitz BSDE
driver satisfying a comparison condition.
A nonlinear domineering claim is a strong
$\mathcal E^f$-martingale dominating the reward
process at every stopping time.
We first establish the existence of such claims
for every optional reward in $\mathcal S^2$.
For left-limited, right-upper semicontinuous rewards,
we prove a local dual representation at every
stopping time: the nonlinear optimal stopping
value equals the essential infimum of the values
of all domineering claims defined on the remaining
stochastic interval.
This infimum is attained by a claim constructed
from the associated reflected BSDE.
We further show that a global domineering claim
attains the value at a prescribed stopping time
if and only if the reflection processes have
accumulated no increase up to that time,
with the appropriate endpoint convention.
When the reward is right-upper semicontinuous
in $\mathcal E^f$-expectation and the family
of contact times is nonempty, its essential
infimum is optimal.
Under strict comparison, we also establish
a necessary contact condition for optimality.
No convexity assumption on the driver is required.
\end{abstract}

\maketitle
\pagestyle{plain}
\thispagestyle{plain}

\section{Introduction}
\label{sec:introduction}

Optimal stopping is a fundamental problem in stochastic
analysis, with important applications in mathematical
finance, stochastic control, and decision-making under
uncertainty. In the classical framework, the value of an
optimal stopping problem is characterized by the Snell
envelope, which is the smallest supermartingale dominating
the reward process. This characterization establishes a
close connection between optimal stopping, martingale
theory, and the Doob--Meyer decomposition.

A nonlinear extension of the classical optimal stopping
problem arises when the conditional expectation is
replaced by a nonlinear evaluation induced by a backward
stochastic differential equation (BSDE). Such evaluations,
known as $f$-expectations, provide a natural framework
for studying optimal stopping under nonlinear dynamics.
In particular, BSDEs driven by Brownian motion and
Poisson random measures allow for the incorporation
of jumps into the underlying stochastic system.
For related developments concerning BSDEs with jumps,
nonlinear evaluations, and dynamic risk measures,
we refer to Quenez and Sulem~\cite{QuenezSulem2013}.

More precisely, let $T>0$ be a fixed time horizon
and let $Y=(Y_t)_{0\leq t\leq T}$ be an optional
reward process. For a stopping time $S$ taking
values in $[0,T]$, we denote by $\mathcal T_{S,T}$
the set of stopping times $\tau$ satisfying
$S\leq\tau\leq T$ almost surely.
The nonlinear optimal stopping value at time $S$
is defined by
\begin{equation}
\label{eq:intro-value}
\mathcal R^f(S)
=
\operatorname*{ess\,sup}_{\tau\in\mathcal T_{S,T}}
\mathcal E^f_{S,\tau}(Y_\tau),
\end{equation}
where $\mathcal E^f_{S,\tau}$ denotes the nonlinear
evaluation induced by a BSDE with driver $f$
between the stopping times $S$ and $\tau$.
When $f=0$, this formulation reduces to the
classical optimal stopping problem under
linear expectation.

A fundamental tool for studying nonlinear optimal
stopping is the theory of reflected backward
stochastic differential equations (RBSDEs).
Unlike standard BSDEs, reflected BSDEs involve
an obstacle and additional increasing processes
that force the solution to remain above the
prescribed obstacle while satisfying suitable
minimality conditions. When the obstacle is
identified with the reward process, the first
component of the reflected BSDE provides a
natural representation of the nonlinear Snell
envelope.

In a Brownian--Poisson filtration, let $W$ be
a Brownian motion and let $\widetilde N$ denote
the compensated Poisson random measure on a
mark space $E$. The reflected BSDE associated
with the driver $f$ and the obstacle $Y$ takes
the form
\begin{equation}
\label{eq:intro-rbsde}
\begin{aligned}
R_t
={}&Y_T
+\int_t^T f(s,R_s,Z_s,k_s)\,ds\\
&-\int_t^T Z_s\,dW_s\\
&-\int_t^T\int_E k_s(e)\,
\widetilde N(ds,de)\\
&+A_T-A_t+B_{T-}-B_{t-},
\qquad 0\leq t\leq T,
\end{aligned}
\end{equation}
subject to the obstacle constraint
\[
R_t\geq Y_t,\qquad 0\leq t\leq T,
\]
and the corresponding minimality conditions.
Here, $Z$ and $k$ are the stochastic integrands,
whereas $A$ and $B$ are nondecreasing reflection
processes, with $A$ predictable and $B$ purely
discontinuous.

Under suitable assumptions on the driver and
the regularity of the reward process, the first
component $R$ of the reflected BSDE coincides
with the nonlinear Snell envelope, namely,
\begin{equation}
\label{eq:intro-snell}
R_S=\mathcal R^f(S)
\end{equation}
for every stopping time $S$.
This connection was established by
Grigorova et al.~\cite{GIOOQ2017} for
left-limited, right-upper semicontinuous
obstacles in a Brownian--Poisson filtration.
The study of nonlinear optimal stopping
was subsequently extended to completely
irregular reward processes in a general
filtration by
Grigorova et al.~\cite{GrigorovaEtAl2019}.
These results show that reflected BSDEs
provide a powerful framework for studying
optimal stopping even when the reward
process is not right-continuous.

An essential feature of the reflected-BSDE
representation is that the nonlinear Snell
envelope is generally a strong
$\mathcal E^f$-supermartingale rather than
a strong $\mathcal E^f$-martingale.
The increasing processes $A$ and $B$
appearing in \eqref{eq:intro-rbsde}
describe the reflection mechanism and
account for the deviation from the
martingale property. Consequently, the
reflected-BSDE decomposition provides
not only a representation of the optimal
stopping value but also important
information about its martingale structure.

On the other hand, an alternative approach
to optimal stopping is based on martingales
that dominate the reward process.
In the classical linear setting,
Jamshidian~\cite{Jamshidian2007}
introduced the notion of domineering claims
and investigated their relationship with
optimal stopping, duality, and the
Doob--Meyer decomposition of the Snell
envelope. This approach suggests a
natural question in the nonlinear framework:
can the reflected-BSDE representation of
the nonlinear Snell envelope be used to
construct a strong $\mathcal E^f$-martingale
that dominates the reward process and
attains the nonlinear optimal stopping value?

The purpose of the present paper is to
investigate this question by establishing
a connection between reflected BSDEs
and nonlinear domineering claims in a
Brownian--Poisson filtration.
We assume that the driver $f$ is Lipschitz
continuous and satisfies a suitable
comparison condition.

In this setting, an
$\mathcal E^f$-domineering claim for the
reward process $Y$ is a strong
$\mathcal E^f$-martingale $D$ satisfying
\[
D_\tau\geq Y_\tau
\qquad
\mathbb P\text{-a.s. for every stopping time }\tau.
\]
Here, the strong $\mathcal E^f$-martingale
property means that
\[
D_\sigma
=
\mathcal E^f_{\sigma,\tau}(D_\tau)
\]
for all stopping times $\sigma\leq\tau$
in the interval on which $D$ is defined.
By the comparison property of
$f$-expectations, every such claim
provides an upper bound for the
corresponding nonlinear optimal
stopping value.

We first establish the existence of
nonlinear domineering claims for
square-integrable optional reward processes,
without imposing pathwise continuity
assumptions. We also investigate the
relationship between optimal stopping
times and the nonlinear Snell envelope,
providing a characterization of optimality
in terms of the martingale property of
the value process.

Furthermore, we study stopping times at
which a domineering claim coincides with
the reward process. Under the assumptions
that the reward is right-upper semicontinuous
in $\mathcal E^f$-expectation and that
the family of contact times is nonempty,
we prove that the essential infimum of
these contact times is an optimal stopping
time. In particular, the associated
domineering claim attains the nonlinear
optimal stopping value.

Our main contribution is a dual
representation of the nonlinear Snell
envelope in terms of domineering claims.
For each stopping time $S$, let
$\mathfrak C_S^f(Y)$ denote the class
of strong $\mathcal E^f$-martingales
defined on the stochastic interval
$[S,T]$ and dominating the reward
process at every stopping time in
this interval.

For a left-limited, right-upper
semicontinuous optional reward process,
we establish the local dual representation
\begin{equation}
\label{eq:intro-duality}
\boxed{
\mathcal R^f(S)
=
R_S
=
\operatorname*{ess\,inf}_{
D\in\mathfrak C_S^f(Y)
}D_S,
}
\end{equation}
for every stopping time $S$.
Moreover, the essential infimum is attained
by an optimal nonlinear domineering claim.

The proof of this result relies
essentially on the reflected-BSDE
characterization of the nonlinear
Snell envelope. More precisely,
starting from the value $R_S$ of
the reflected BSDE at time $S$,
we construct a strong
$\mathcal E^f$-martingale
$D^{*,S}$ using the stochastic
integrands $(Z,k)$ appearing in
\eqref{eq:intro-rbsde}.
The construction is designed so that
\[
D_S^{*,S}=R_S.
\]
An integrating-factor argument involving
the reflection processes then yields
\[
D_t^{*,S}\geq R_t\geq Y_t,
\qquad S\leq t\leq T,
\]
in the appropriate almost-sure sense.
Consequently, $D^{*,S}$ is an
optimal domineering claim and attains
the dual representation
\eqref{eq:intro-duality}.
This construction does not require
any convexity assumption on the driver.

We further investigate the distinction
between local and global dual attainment.
While the local dual representation
admits an optimal domineering claim
at every stopping time, the attainment
by a claim defined on the entire
interval $[0,T]$ is closely related
to the reflection mechanism of the
associated RBSDE.

To make this connection precise,
let $\mathfrak C^f(Y)$ denote the
global class of strong
$\mathcal E^f$-martingales on $[0,T]$
that dominate the reward process at
every stopping time. We establish
the following characterization:
for every stopping time $S$,
\begin{equation}
\label{eq:intro-global-duality}
\boxed{
\begin{aligned}
&\exists D\in\mathfrak C^f(Y)
\text{ such that }D_S=R_S\\
&\quad\Longleftrightarrow
R\text{ is a strong }
\mathcal E^f\text{-martingale on }[0,S]\\
&\quad\Longleftrightarrow
A_S=B_{S-}=0
\quad\mathbb P\text{-a.s.}
\end{aligned}
}
\end{equation}

This equivalence shows that global
dual attainment at a prescribed stopping
time is determined by the absence of
reflection up to that time, in the
precise sense of
\eqref{eq:intro-global-duality}.
It therefore highlights the fundamental
role of the increasing processes in the
reflected-BSDE decomposition and reveals
a direct relationship between the
martingale structure of the nonlinear
Snell envelope and the attainment of
the global dual problem.

Finally, we establish a necessary
contact condition for optimality.
Under a suitable strict comparison
assumption, we show that if a global
domineering claim attains the nonlinear
optimal stopping value at the initial
time, then it must coincide with the
reward process at any optimal stopping
time. This result provides an additional
connection between dual attainment,
optimal stopping, and the contact
structure of the reward process.

The paper is organized as follows.
Section~2 introduces the probabilistic
framework, BSDEs with jumps, and
$f$-expectations. Section~3 recalls
the reflected-BSDE characterization
of the nonlinear Snell envelope
and presents a characterization of
optimal stopping times. Section~4
introduces nonlinear domineering
claims, establishes their existence,
and investigates contact times.
It also develops the local dual
representation, characterizes global
dual attainment through the reflection
processes, and establishes a necessary
contact condition for optimality.

\section{Probabilistic framework and \texorpdfstring{$f$}{f}-expectation}
\label{sec:framework}
Let $T>0$, and let $(E,\mathcal E)$ be a measurable
space equipped with a $\sigma$-finite measure $\nu$.
Let $(\Omega,\mathcal F,\mathbb P)$ be a complete
probability space supporting a one-dimensional
Brownian motion $W=(W_t)_{t\in[0,T]}$ and an
independent Poisson random measure
$N(dt,de)$ on $[0,T]\times E$, with compensator
$dt\,\nu(de)$.

We assume that $\mathbb F=(\mathcal F_t)_{t\in[0,T]}$
is the usual augmentation of the natural filtration
generated by $W$ and $N$.
In particular, $\mathbb F$ satisfies the usual
conditions and enjoys the predictable representation
property.

The compensated Poisson random measure is defined by
\[
\widetilde N(dt,de)
:=
N(dt,de)-dt\,\nu(de).
\]

Let $\mathcal{T}_0$ denote the set of all $\mathbb{F}$-stopping times
with values in $[0,T]$. For each $t\in[0,T]$, we set
\[
\mathcal{T}_{t,T}
:=
\left\{
\tau\in\mathcal{T}_0:
t\leq \tau\leq T
\quad \mathbb{P}\text{-a.s.}
\right\}.
\]

For every $S\in\mathcal{T}_0$, we define
\[
\mathcal{T}_S
:=
\left\{
\theta\in\mathcal{T}_0:
\theta\geq S
\quad \mathbb{P}\text{-a.s.}
\right\}.
\]

We shall use the following notation and functional spaces.

\begin{itemize}

\item $\mathcal{P}$ denotes the predictable $\sigma$-algebra on
$\Omega\times[0,T]$.

\item $L^{2}(\mathcal{F}_T)$ denotes the space of all
$\mathcal{F}_T$-measurable random variables $\xi$ such that
\[
\|\xi\|_{L^{2}(\mathcal{F}_T)}^{2}
:=
\mathbb{E}\bigl[|\xi|^{2}\bigr]
<\infty.
\]

\item $L_{\nu}^{2}$ denotes the Hilbert space
\[
L_{\nu}^{2}
:=
\left\{
\ell:E\longrightarrow\mathbb{R}
\text{ measurable}:
\|\ell\|_{\nu}^{2}
:=
\int_E |\ell(e)|^{2}\,\nu(\mathrm{d}e)
<\infty
\right\}.
\]
It is equipped with the inner product $\langle \ell,\kappa\rangle_{\nu}
:=
\int_E \ell(e)\kappa(e)\,\nu(\mathrm{d}e).
$

\item $\mathcal{B}(L_{\nu}^{2})$ denotes the Borel
$\sigma$-algebra on $L_{\nu}^{2}$.

\item $\mathbb{H}^{2}$ denotes the space of all real-valued
$\mathbb{F}$-predictable processes $X$ such that
\[
\|X\|_{\mathbb{H}^{2}}^{2}
:=
\mathbb{E}\left[
\int_0^T |X_t|^{2}\,\mathrm{d}t
\right]
<\infty.
\]

\item $\mathbb{H}_{\nu}^{2}$ denotes the space of all
$\mathcal{P}\otimes\mathcal{E}$-measurable processes
\[
\ell:\Omega\times[0,T]\times E\longrightarrow\mathbb{R}
\]
such that
\[
\|\ell\|_{\mathbb{H}_{\nu}^{2}}^{2}
:=
\mathbb{E}\left[
\int_0^T \|\ell_t\|_{\nu}^{2}\,\mathrm{d}t
\right]
=
\mathbb{E}\left[
\int_0^T\int_E
|\ell_t(e)|^{2}\,\nu(\mathrm{d}e)\,\mathrm{d}t
\right]
<\infty.
\]

\item $\mathcal{S}^{2}$ denotes the space of all real-valued optional
processes $X=(X_t)_{t\in[0,T]}$ such that
\[
\|X\|_{\mathcal{S}^{2}}^{2}
:=
\mathbb{E}\left[
\operatorname*{ess\,sup}_{\tau\in\mathcal{T}_0}
|X_\tau|^{2}
\right]
<\infty.
\]

\item $\mathcal{M}^{2}$ denotes the Hilbert space of all
square-integrable martingales
$M=(M_t)_{t\in[0,T]}$ satisfying $M_0=0$. It is equipped with the
inner product
\[
\langle M,M'\rangle_{\mathcal{M}^{2}}
:=
\mathbb{E}[M_T M'_T]
=
\mathbb{E}\bigl[\langle M,M'\rangle_T\bigr]
=
\mathbb{E}\bigl[[M,M']_T\bigr],
\]
and the norm $\|M\|_{\mathcal{M}^{2}}^{2}
:=
\mathbb{E}\bigl[|M_T|^{2}\bigr]$.

\end{itemize}

\begin{definition}[Driver and Lipschitz driver]
\label{def:driver}

A function
\[
g:\Omega\times[0,T]\times\mathbb{R}^{2}\times L_{\nu}^{2}
\longrightarrow \mathbb{R},
\qquad
(\omega,t,y,z,k)
\longmapsto
g(\omega,t,y,z,k),
\]
is called a \emph{driver} if it satisfies the following conditions:

\begin{enumerate}
\item[\textnormal{(i)}] \textbf{Measurability:}
the function $g$ is $\mathcal{P}
\otimes\mathcal{B}(\mathbb{R}^{2})
\otimes\mathcal{B}(L_{\nu}^{2})
\text{-measurable}$.

\item[\textnormal{(ii)}] \textbf{Integrability:}
\[
\mathbb{E}\left[
\int_{0}^{T}
|g(t,0,0,0)|^{2}\,\mathrm{d}t
\right]
<\infty.
\]
\end{enumerate}

A driver $g$ is called a \emph{Lipschitz driver} if there exists a
constant $K\geq 0$ such that, $\mathrm{d}\mathbb{P}\otimes\mathrm{d}t$-a.e.,
for every
\[
(y_i,z_i,k_i)\in\mathbb{R}^{2}\times L_{\nu}^{2},
\qquad i\in\{1,2\},
\]
we have
\[
\begin{aligned}
&\left|
g(\omega,t,y_1,z_1,k_1)
-
g(\omega,t,y_2,z_2,k_2)
\right|
\\
&\leq
K\left(
|y_1-y_2|
+
|z_1-z_2|
+
\|k_1-k_2\|_{\nu}
\right).
\end{aligned}
\]

A pair $(g,\xi)$, where $g$ is a Lipschitz driver and
$\xi\in L^{2}(\mathcal{F}_T)$, is called \emph{standard data}.
\end{definition}

\begin{definition}[BSDE and conditional \texorpdfstring{$g$}{g}-expectation]
\label{def:bsde-g-expectation}

Let $(g,\xi)$ be standard data. The backward stochastic differential
equation with driver $g$, terminal time $T$, and terminal condition
$\xi$ is given by
\begin{equation}
\label{eq:bsde-g}
\begin{aligned}
X_t
={}&
\xi
+
\int_t^T g(s,X_s,\pi_s,l_s)\,\mathrm{d}s
-
\int_t^T \pi_s\,\mathrm{d}W_s
\\
&-
\int_t^T\int_E
l_s(e)\,
\widetilde{N}(\mathrm{d}s,\mathrm{d}e)
,
\qquad 0\leq t\leq T.
\end{aligned}
\end{equation}

Under the Lipschitz assumption on $g$, \eqref{eq:bsde-g} admits a
unique solution
\[
(X,\pi,l)
\in
\mathcal{S}^{2}
\times\mathbb{H}^{2}
\times\mathbb{H}_{\nu}^{2}.
\]

For every $t\in[0,T]$, the nonlinear operator
\[
\mathcal{E}_{t,T}^{g}:
L^{2}(\mathcal{F}_T)
\longrightarrow
L^{2}(\mathcal{F}_t)
\]
defined by
\[
\mathcal{E}_{t,T}^{g}(\xi):=X_t
\]
is called the \emph{conditional $g$-expectation} of $\xi$ at time $t$.

This definition extends to random horizons. More precisely, if
$\sigma$ and $\tau$ are stopping times such that
\[
0\leq \sigma\leq\tau\leq T
\qquad \mathbb{P}\text{-a.s.},
\]
then
\[
\mathcal{E}_{\sigma,\tau}^{g}:
L^{2}(\mathcal{F}_{\tau})
\longrightarrow
L^{2}(\mathcal{F}_{\sigma})
\]
is defined by
\[
\mathcal{E}_{\sigma,\tau}^{g}(\eta)
:=
X_{\sigma},
\]
where $(X,\pi,l)$ is the unique solution, on the stochastic interval
$[\![\sigma,\tau]\!]$, of the BSDE with terminal condition
$\eta\in L^{2}(\mathcal{F}_{\tau})$.
\end{definition}

We impose the following assumption on the driver $g$. This condition
ensures the comparison property for BSDEs with jumps; see
Quenez and Sulem~\cite[Theorem~4.2]{QuenezSulem2013}.

\begin{assumptionA}
\label{ass:comparison}

There exists a mapping
\[
\gamma:
[0,T]\times\Omega\times
\mathbb{R}^{2}\times
(L_{\nu}^{2})^{2}
\longrightarrow
L_{\nu}^{2},
\]
written as
\[
(t,\omega,y,z,k_1,k_2)
\longmapsto
\gamma_t^{y,z,k_1,k_2}(\omega,\cdot),
\]
which is
\[
\mathcal{P}
\otimes\mathcal{B}(\mathbb{R}^{2})
\otimes\mathcal{B}(L_{\nu}^{2})^{\otimes 2}
\text{-measurable},
\]
and such that, $\mathrm{d}\mathbb{P}\otimes\mathrm{d}t$-a.e., for every
\[
(y,z,k_1,k_2)
\in
\mathbb{R}^{2}\times(L_{\nu}^{2})^{2},
\]
we have
\begin{equation}
\label{eq:assumption-A}
g(t,y,z,k_1)-g(t,y,z,k_2)
\geq
\left\langle
\gamma_t^{y,z,k_1,k_2},
k_1-k_2
\right\rangle_{\nu}.
\end{equation}

Moreover, there exists a constant $K>0$ such that
\[
\left\|
\gamma_t^{y,z,k_1,k_2}
\right\|_{\nu}
\leq K,
\qquad
\mathrm{d}\mathbb{P}\otimes\mathrm{d}t\text{-a.e.},
\]
and
\[
\gamma_t^{y,z,k_1,k_2}(e)\geq -1,
\qquad
\mathrm{d}\mathbb{P}\otimes\mathrm{d}t
\otimes\nu(\mathrm{d}e)\text{-a.e.}
\]
\end{assumptionA}
\section{The nonlinear Snell envelope}
\label{sec:snell}

From now on, $f$ denotes a predictable Lipschitz driver satisfying
Assumption~A.

Let $Y=(Y_t)_{0\leq t\leq T}\in\Sspace$ be an optional reward process. For
each $S\in\T_{0,T}$, define the value family
\begin{equation}
\label{eq:value-family}
\Snell(S):=\esssup_{\tau\in\T_{S}}\Ef{S}{\tau}(Y_\tau).
\end{equation}
\subsection{Reflected BSDE and Snell-envelope characterization}
To avoid confusion with the domineering claim introduced later, the purely
discontinuous increasing process in the reflected BSDE is denoted by $B$.

\begin{definition}[Reflected BSDE]
\label{def:rbsde}
Let $Y\in\mathcal S^2$ be an obstacle. A quintuple
$(R,Z,k,A,B)$ is called a solution of the reflected BSDE with parameters
$(f,Y)$ if
\[
(R,Z,k,A,B)
\in
\mathcal S^2\times\mathbb H^2\times\mathbb H_\nu^2
\times\mathcal S^2\times\mathcal S^2,
\]
and, almost surely, for every $t\in[0,T]$,
\begin{equation}
\label{eq:rbsde}
\begin{aligned}
R_t={}&Y_T+\int_t^T f(s,R_s,Z_s,k_s)\,\mathrm ds
-\int_t^T Z_s\,\mathrm dW_s\\
&-\int_t^T\int_E k_s(e)\,\widetilde N(\mathrm ds,\mathrm de)
+A_T-A_t+B_{T-}-B_{t-},
\end{aligned}
\end{equation}
with
\[
R_t\geq Y_t,\qquad 0\leq t\leq T.
\]
The process $A$ is nondecreasing, right-continuous and predictable, with
$A_0=0$, and satisfies
\[
\int_0^T\ind_{\{R_t>Y_t\}}\,\mathrm dA_t^c=0,
\qquad
(R_{\tau-}-Y_{\tau-})\Delta A_\tau^d=0
\]
for every predictable stopping time $\tau$. The process $B$ is nondecreasing,
right-continuous, adapted and purely discontinuous, with $B_{0-}=0$, and
satisfies
\[
(R_\tau-Y_\tau)\Delta B_\tau=0
\]
for every $\tau\in\mathcal T_{0,T}$.
\end{definition}
The following result is \cite[Theorem~4.1]{GIOOQ2017}, rewritten with the present notation.
\begin{theorem}[Reflected-BSDE characterization]
\label{thm:rbsde-characterization}
Assume that $Y\in\mathcal S^2$ is left-limited and right-upper
semicontinuous. Let $(R,Z,k,A,B)$ be the solution of the reflected BSDE with
parameters $(f,Y)$. Then, for every $S\in\mathcal T_{0,T}$,
\begin{equation}
\label{eq:rbsde-value}
R_S
=\operatorname*{ess\,sup}_{\tau\in\mathcal T_{S,T}}
\mathcal E^f_{S,\tau}(Y_\tau)
=\mathcal R^f(S)
\qquad\mathbb P\text{-a.s.}
\end{equation}
Moreover, for every $\varepsilon>0$, the stopping time
\[
\tau_S^\varepsilon
:=\inf\{t\geq S:R_t\leq Y_t+\varepsilon\}\wedge T
\]
is $L\varepsilon$-optimal: there exists a constant $L$, depending only on
$T$ and the Lipschitz constant of $f$, such that
\begin{equation}
\label{eq:epsilon-optimality}
R_S
\leq
\mathcal E^f_{S,\tau_S^\varepsilon}
\bigl(Y_{\tau_S^\varepsilon}\bigr)
+L\varepsilon
\qquad\mathbb P\text{-a.s.}
\end{equation}
\end{theorem}
\subsection{Optimal stopping times}

\begin{definition}[Optimal stopping time]
A stopping time $\tau^{*}\in\T_{S,T}$ is called $S$-optimal if
\[
\Snell(S)=\Ef{S}{\tau^{*}}(Y_{\tau^{*}})\qquad\Pp\text{-a.s.}
\]
\end{definition}
The argument used in the following proposition is inspired by
Bouhadou and Ouknine~\cite[Section~4, Proposition~.2]{BouhadouOuknine2021},
who study a related characterization for strong supermartingales
under predictable nonlinear expectations.
Here, we adapt this idea to the optional setting with
$f$-expectations in a Brownian--Poisson filtration.

\begin{proposition}[Characterization of optimal stopping times]
\label{prop:optimality-characterization}

Let $Y\in\mathcal S^2$ be an optional process,
and let $S\in\mathcal T_{0,T}$ and
$\tau^*\in\mathcal T_S$.

Then the following assertions are equivalent:

\begin{enumerate}[label=\textnormal{(\arabic*)}]

\item $\tau^*$ is $S$-optimal, that is,
\[
\mathcal R^f(S)
=
\mathcal E^f_{S,\tau^*}(Y_{\tau^*})
\quad\mathbb P\text{-a.s.}
\]

\item We have
\[
\mathcal R^f(S)
=
\mathcal E^f_{S,\tau^*}
\bigl(\mathcal R^f(\tau^*)\bigr)
=
\mathcal E^f_{S,\tau^*}(Y_{\tau^*})
\quad\mathbb P\text{-a.s.}
\]

\end{enumerate}

If, in addition, the evaluation
$\mathcal E^f_{S,\tau^*}$ satisfies
the strict comparison property, these
assertions are also equivalent to:

\begin{enumerate}[label=\textnormal{(\arabic*)},
start=3]
\item
\[
\mathcal R^f(\tau^*)=Y_{\tau^*}
\quad\mathbb P\text{-a.s.},
\]
and
\[
\mathcal R^f(S)
=
\mathcal E^f_{S,\tau^*}
\bigl(\mathcal R^f(\tau^*)\bigr)
\quad\mathbb P\text{-a.s.}
\]
\end{enumerate}

\end{proposition}

\begin{proof}
Since $\mathcal R^f$ is a strong
$\mathcal E^f$-supermartingale family,
and $\mathcal R^f(\tau^*)\geq Y_{\tau^*}$,
the comparison property implies
\[
\mathcal R^f(S)
\geq
\mathcal E^f_{S,\tau^*}
\bigl(\mathcal R^f(\tau^*)\bigr)
\geq
\mathcal E^f_{S,\tau^*}(Y_{\tau^*}).
\]

Consequently, assertion (1) holds
if and only if both inequalities
are equalities, which is precisely
assertion (2).

Suppose now that the strict comparison
property holds for
$\mathcal E^f_{S,\tau^*}$.

If (2) holds, then
\[
\mathcal E^f_{S,\tau^*}
\bigl(\mathcal R^f(\tau^*)\bigr)
=
\mathcal E^f_{S,\tau^*}(Y_{\tau^*}).
\]
Since
$\mathcal R^f(\tau^*)\geq Y_{\tau^*}$,
strict comparison yields
\[
\mathcal R^f(\tau^*)=Y_{\tau^*}
\quad\mathbb P\text{-a.s.}
\]
Thus (3) follows.

Conversely, if (3) holds, then
\[
\mathcal R^f(S)
=
\mathcal E^f_{S,\tau^*}
\bigl(\mathcal R^f(\tau^*)\bigr)
=
\mathcal E^f_{S,\tau^*}(Y_{\tau^*}),
\]
so $\tau^*$ is $S$-optimal.
\end{proof}

\section{Domineering claims and the duality result}
\begin{definition}[Nonlinear domineering claim]
\label{def:domineering}
A process $C\in\Sspace$ is called an $\mathcal{E}^{f}$-domineering claim for
$Y$ if
\begin{enumerate}[label=(\roman*)]
\item $C$ is a strong $\mathcal{E}^{f}$-martingale, that is,
\[
C_\sigma=\Ef{\sigma}{\tau}(C_\tau)
\quad\text{for all }\sigma,\tau\in\T_{0,T},\ \sigma\leq\tau;
\]
\item $C_\tau\geq Y_\tau$ almost surely for every
$\tau\in\T_{0,T}$.
\end{enumerate}
\end{definition}
\begin{theorem}[Existence of an $\mathcal E^f$-domineering claim]
\label{thm:existence-domineering}
Let $Y\in\mathcal S^2$ be an optional process. Then there exists an
$\mathcal E^f$-domineering claim $C$ for $Y$.

More precisely, there exists $H\in L^2(\mathcal F_T)$ such that
\begin{equation}
\label{eq:domineering-construction}
C_t:=\mathcal E^f_{t,T}(H),\qquad 0\leq t\leq T,
\end{equation}
is a strong $\mathcal E^f$-martingale and
\[
C_\tau\geq Y_\tau
\qquad\mathbb P\text{-a.s. for every }\tau\in\mathcal T_{0,T}.
\]

\end{theorem}

\begin{proof}
Let $Y\in\mathcal S^2$ be an optional process.

\medskip
\noindent
\textbf{Step 1: Construction of terminal variables.}

Fix $\tau\in\mathcal T_{0,T}$.
Consider the random ordinary differential equation
on $[\![\tau,T]\!]$:
\[
X_s^\tau
=
Y_\tau-\int_\tau^s
f(r,X_r^\tau,0,0)\,dr,
\qquad \tau\leq s\leq T.
\]

Since $f$ is predictable and uniformly Lipschitz
in its value variable, this equation admits a
unique adapted, continuous solution on
$[\![\tau,T]\!]$.

Set
\[
\eta^\tau:=X_T^\tau.
\]

Then $\eta^\tau$ is $\mathcal F_T$-measurable.
Rewriting the equation backward, we obtain
\[
X_s^\tau
=
\eta^\tau+
\int_s^T f(r,X_r^\tau,0,0)\,dr,
\qquad \tau\leq s\leq T.
\]

Thus $(X^\tau,0,0)$ solves the BSDE on
$[\![\tau,T]\!]$ with terminal condition
$\eta^\tau$. In particular, once the
square-integrability established below is used,
uniqueness yields
\begin{equation}
\label{eq:terminal-identity}
\mathcal E^f_{\tau,T}(\eta^\tau)
=
X_\tau^\tau
=
Y_\tau.
\end{equation}

\medskip
\noindent
\textbf{Step 2: Uniform square-integrability.}

Define
\[
Y^*:=
\operatorname*{ess\,sup}_{\theta\in\mathcal T_{0,T}}
|Y_\theta|
\]
and
\[
F^*:=\int_0^T|f(r,0,0,0)|\,dr.
\]

By the assumptions on $Y$ and $f$,
\[
Y^*\in L^2(\mathcal F_T),
\qquad
F^*\in L^2(\mathcal F_T).
\]

Indeed, by Cauchy--Schwarz,
\[
\mathbb E[(F^*)^2]
\leq
T\mathbb E\left[
\int_0^T|f(r,0,0,0)|^2\,dr
\right]
<\infty.
\]

Using the Lipschitz property,
\[
|X_s^\tau|
\leq
|Y_\tau|
+
\int_\tau^s |f(r,0,0,0)|\,dr
+
K\int_\tau^s |X_r^\tau|\,dr.
\]

Gronwall's inequality gives
\[
|X_s^\tau|
\leq e^{KT}(Y^*+F^*),
\qquad \tau\leq s\leq T.
\]

In particular,
\begin{equation}
\label{eq:terminal-uniform-bound}
|\eta^\tau|
\leq e^{KT}(Y^*+F^*)
\quad\mathbb P\text{-a.s.}
\end{equation}

This estimate also ensures the required
square-integrability of the solution
$X^\tau$ on its stochastic interval.

\medskip
\noindent
\textbf{Step 3: Construction of a common terminal condition.}

Define
\[
H:=
\operatorname*{ess\,sup}_{\tau\in\mathcal T_{0,T}}
\eta^\tau.
\]

The essential supremum is well defined
as an $\mathcal F_T$-measurable random variable.
By \eqref{eq:terminal-uniform-bound}

\[
-e^{KT}(Y^*+F^*)
\leq H
\leq e^{KT}(Y^*+F^*).
\]

Consequently,
\[
H\in L^2(\mathcal F_T).
\]

\medskip
\noindent
\textbf{Step 4: Construction of the domineering claim.}

Let $(C,Z^C,k^C)$ be the unique solution
of the BSDE with driver $f$ and terminal
condition $H$.

Then
\[
C_t=\mathcal E^f_{t,T}(H),
\qquad 0\leq t\leq T,
\]
and $C\in\mathcal S^2$.

For any stopping times
$\sigma,\tau\in\mathcal T_{0,T}$
such that $\sigma\leq\tau$, the time
consistency property yields
\[
\begin{aligned}
\mathcal E^f_{\sigma,\tau}(C_\tau)
&=
\mathcal E^f_{\sigma,\tau}
\left(\mathcal E^f_{\tau,T}(H)\right)\\
&=
\mathcal E^f_{\sigma,T}(H)\\
&=C_\sigma.
\end{aligned}
\]

Hence $C$ is a strong
$\mathcal E^f$-martingale.

Moreover, for every
$\tau\in\mathcal T_{0,T}$, the definition
of $H$ implies
\[
H\geq\eta^\tau
\quad\mathbb P\text{-a.s.}
\]

By the comparison theorem for BSDEs with
jumps and \eqref{eq:terminal-identity},
\[
\begin{aligned}
C_\tau
&=\mathcal E^f_{\tau,T}(H)\\
&\geq\mathcal E^f_{\tau,T}(\eta^\tau)\\
&=Y_\tau
\quad\mathbb P\text{-a.s.}
\end{aligned}
\]

Therefore, $C$ is an
$\mathcal E^f$-domineering claim for $Y$.
\end{proof}

\subsection{Dual upper bound and contact time}

We now establish a nonlinear dual upper bound. We then show that this bound is
attained whenever the family of contact times between the payoff process and
the domineering claim is nonempty.\\
We introduce the following regularity property, which will be used in the
study of optimal stopping times associated with an $\mathcal E^f$-domineering
claim.

\begin{definition}[Right-upper semicontinuity in
$\mathcal E^f$-expectation]
\label{def:RUSC-Ef}
An optional process $Y\in\Sspace$ is said to be
\emph{right-upper semicontinuous in $\mathcal E^f$-expectation}
if, for every decreasing sequence
$(\tau_n)_{n\geq1}\subset\T_{0,T}$ such that
$\tau_n\downarrow\tau$ $\Pp$-a.s. for some
$\tau\in\T_{0,T}$, we have
\[
\Ef{0}{\tau}(Y_\tau)
\geq
\limsup_{n\to\infty}
\Ef{0}{\tau_n}(Y_{\tau_n})
\qquad \Pp\text{-a.s.}
\]
\end{definition}
Motivated by the classical duality approach of Jamshidian \cite{Jamshidian2007}, we investigate the relationship between nonlinear domineering claims and optimal stopping times in the framework of \(\mathcal E^f\)-expectations. Our aim is to extend the contact-time characterization of optimal stopping to optional payoff processes satisfying right-upper semicontinuity in \(\mathcal E^f\)-expectation, without requiring pathwise right-continuity. In particular, the following proposition shows that, whenever a strong \(\mathcal E^f\)-martingale dominates the payoff and the corresponding contact family is nonempty, the minimal contact time is optimal and the nonlinear duality equality holds.

\begin{proposition}[Optimal stopping induced by a domineering claim]
\label{prop:duality-contact-time}

Let $Y\in\mathcal S^2$ be an optional payoff process, $C\in\mathcal S^2$ be an
$\mathcal E^f$-domineering claim for $Y$ and $S \in \mathcal T_{0,T}$.

Then
\begin{equation}
\label{eq:dual-upper-bound}
\operatorname*{ess\,sup}_{\tau\in\mathcal T_{S,T}}
\mathcal E^f_{S,\tau}(Y_\tau)
\leq C_S.
\end{equation}

Suppose, in addition, that $Y$ is
right-upper semicontinuous in
$\mathcal E^f$-expectation and that
\[
\mathscr A_C
:=
\left\{
\sigma\in\mathcal T_{0,T}:
C_\sigma=Y_\sigma\quad\mathbb P\text{-a.s.}
\right\}
\]
is nonempty.

Define
\[
T_C^*
:=
\operatorname*{ess\,inf}_{\sigma\in\mathscr A_C}
\sigma.
\]

Then $T_C^*\in\mathcal T_{0,T}$ is an optimal
stopping time and
\begin{equation}
\label{eq:duality-equality-C}
\operatorname*{ess\,sup}_{\tau\in\mathcal T_{0,T}}
\mathcal E^f_{0,\tau}(Y_\tau)
=
\mathcal E^f_{0,T_C^*}(Y_{T_C^*})
=
C_0.
\end{equation}
\end{proposition}

\begin{example}[A simple case where $\mathscr{A}_C\neq\varnothing$]
\label{ex:nonempty-contact-set}

Let $f=0$. Then the $f$-expectation reduces to the classical conditional
expectation. Fix $T>0$ and consider
\[
C_t:=1,
\qquad
Y_t:=\frac{t}{T},
\qquad 0\leq t\leq T.
\]

The process $C$ is a martingale and
\[
C_t=1\geq \frac{t}{T}=Y_t,
\qquad 0\leq t\leq T.
\]
Hence $C$ is an $\mathcal{E}^{0}$-domineering claim for $Y$.

Moreover,
\[
C_\sigma=Y_\sigma
\quad\Longleftrightarrow\quad
1=\frac{\sigma}{T}
\quad\Longleftrightarrow\quad
\sigma=T
\qquad \Pp\text{-a.s.}
\]
Therefore,
\[
\mathscr{A}_C=\{T\},
\]
so that, in particular,
\[
\mathscr{A}_C\neq\varnothing.
\]

Furthermore,
\[
T_C^*
=
\essinf_{\sigma\in\mathscr{A}_C}\sigma
=
T.
\]
Finally,
\[
\esssup_{\tau\in\T_{0,T}}\mathbb{E}[Y_\tau]
=
\esssup_{\tau\in\T_{0,T}}
\mathbb{E}\left[\frac{\tau}{T}\right]
=
1
=
C_0,
\]
and hence the optimal stopping time is $T_C^*=T$.
\end{example}

\begin{proof}
For any $\tau\in\mathcal T_{S,T}$,
the domination property and the comparison
theorem yield
\[
\mathcal E^f_{S,\tau}(Y_\tau)
\leq
\mathcal E^f_{S,\tau}(C_\tau)
=C_S.
\]

Taking the essential supremum gives
\begin{equation}
\label{eq:contact-upper-bound}
\operatorname*{ess\,sup}_{\tau\in\mathcal T_{S,T}}
\mathcal E^f_{S,\tau}(Y_\tau)
\leq C_S.
\end{equation}

We next show that $T_C^*$ is a stopping time.

Let $\sigma_1,\sigma_2\in\mathscr A_C$.
Since
\[
C_{\sigma_i}=Y_{\sigma_i}
\quad\mathbb P\text{-a.s.},
\qquad i=1,2,
\]
we obtain
\[
\begin{aligned}
C_{\sigma_1\wedge\sigma_2}
&=
C_{\sigma_1}
\mathbf 1_{\{\sigma_1\leq\sigma_2\}}
+
C_{\sigma_2}
\mathbf 1_{\{\sigma_2<\sigma_1\}}\\
&=
Y_{\sigma_1\wedge\sigma_2}
\quad\mathbb P\text{-a.s.}
\end{aligned}
\]

Thus $\mathscr A_C$ is stable under
pairwise minima.

Consequently, the family $\mathscr A_C$
is directed downwards. By the standard
essential-infimum approximation theorem,
there exists a decreasing sequence
$(\sigma_n)_{n\geq1}\subset\mathscr A_C$
such that
\[
\sigma_n\downarrow T_C^*
\quad\mathbb P\text{-a.s.}
\]

Since the filtration is right-continuous,
$T_C^*$ is an $\mathbb F$-stopping time.
Moreover, $0\leq T_C^*\leq T$, and hence
\[
T_C^*\in\mathcal T_{0,T}.
\]

For each $n\geq1$, since
$\sigma_n\in\mathscr A_C$, we have
\[
Y_{\sigma_n}=C_{\sigma_n}
\quad\mathbb P\text{-a.s.}
\]

Therefore, the strong
$\mathcal E^f$-martingale property yields
\[
\mathcal E^f_{0,\sigma_n}(Y_{\sigma_n})
=
\mathcal E^f_{0,\sigma_n}(C_{\sigma_n})
=C_0.
\]

By the right-upper semicontinuity
of $Y$ in $\mathcal E^f$-expectation,
\[
\begin{aligned}
\mathcal E^f_{0,T_C^*}(Y_{T_C^*})
&\geq
\limsup_{n\to\infty}
\mathcal E^f_{0,\sigma_n}(Y_{\sigma_n})\\
&=C_0.
\end{aligned}
\]

On the other hand, since $C$ dominates
$Y$, the comparison theorem and
the strong martingale property give
\[
\begin{aligned}
\mathcal E^f_{0,T_C^*}(Y_{T_C^*})
&\leq
\mathcal E^f_{0,T_C^*}(C_{T_C^*})\\
&=C_0.
\end{aligned}
\]

It follows that
\[
\mathcal E^f_{0,T_C^*}(Y_{T_C^*})
=C_0.
\]

Combining this equality with
\eqref{eq:contact-upper-bound} for $S=0$, we obtain
\[
\boxed{
\operatorname*{ess\,sup}_{\tau\in\mathcal T_{0,T}}
\mathcal E^f_{0,\tau}(Y_\tau)
=
\mathcal E^f_{0,T_C^*}(Y_{T_C^*})
=C_0.
}
\]

Hence $T_C^*$ is optimal.
\end{proof}

\begin{remark}
The optimality of $T_C^*$ does not
necessarily imply that
\[
C_{T_C^*}=Y_{T_C^*}
\quad\mathbb P\text{-a.s.}
\]
under the comparison assumption alone.

However, if the evaluation
$\mathcal E^f_{0,T_C^*}$ satisfies
the strict comparison property, then
\[
\mathcal E^f_{0,T_C^*}(C_{T_C^*})
=
\mathcal E^f_{0,T_C^*}(Y_{T_C^*})
\]
together with
$C_{T_C^*}\geq Y_{T_C^*}$
implies
\[
C_{T_C^*}=Y_{T_C^*}
\quad\mathbb P\text{-a.s.}
\]

Thus, under strict comparison,
$T_C^*$ belongs to $\mathscr A_C$
and is its smallest element
in the almost-sure order.
\end{remark}

We now establish a nonlinear dual representation
of the optimal stopping value, extending the
classical domineering-claim approach of
Jamshidian \cite{Jamshidian2007} to the
framework of $\mathcal E^f$-expectations.

Unlike the preceding result, which characterizes
optimal stopping in terms of a given domineering
claim, the following proposition addresses the
existence of an optimal domineering claim.

More precisely, for a left-limited and
right-upper semicontinuous optional payoff
process, we show that the nonlinear optimal
stopping value coincides with the essential
infimum of the initial values of all strong
$\mathcal E^f$-martingales dominating the payoff.

Furthermore, this infimum is attained, yielding
a nonlinear analogue of the classical duality
principle in a framework driven by Brownian
motion and Poisson jumps.

\begin{proposition}[Dual representation]
\label{prop:dual-representation}

Let \(Y\in\mathcal{S}^{2}\) be a left-limited and right-upper
semicontinuous optional process. Let \(f\) be a Lipschitz driver
satisfying the comparison assumption, and let \(S\) be a stopping
time in \(\mathcal{T}_{0,T}\).
Let \(\mathfrak{C}_S^{f}(Y)\) be the set of all strong
\(\mathcal{E}^{f}\)-martingales \(D\in\mathcal{S}^{2}\) on
\([S,T]\) such that
\[
D_\tau\geq Y_\tau
\quad\text{a.s. for every }\tau\in\mathcal{T}_{S,T}.
\]

Then \(\mathfrak{C}_S^{f}(Y)\neq\varnothing\), and
\[
\boxed{
\operatorname*{ess\,sup}_{\tau\in\mathcal{T}_{S,T}}
\mathcal{E}^{f}_{S,\tau}(Y_\tau)
=
R_S
=
\operatorname*{ess\,inf}_{D\in\mathfrak{C}_S^{f}(Y)}D_S.
}
\]
Moreover, the essential infimum is attained by a claim \(D^*\).
\end{proposition}

\begin{proof}
We first identify the value process by applying
Theorem~3.2 to the payoff process $Y$.
Let $(R,Z,k,A,B)$ be the solution of the reflected
BSDE associated with $(f,Y)$.
Then, for every $S\in\mathcal T_{0,T}$,
\[
R_S
=
\mathcal R^f(S)
=
\operatorname*{ess\,sup}_{\tau\in\mathcal T_{S,T}}
\mathcal E^f_{S,\tau}(Y_\tau).
\]

Moreover, by the obstacle condition,
\[
R_\tau\geq Y_\tau,
\qquad
\text{a.s. for every }\tau\in\mathcal T_{0,T}.
\]

Let $D\in\mathfrak C_S^f(Y)$.
For every $\tau\in\mathcal T_{S,T}$,
the domination property and the comparison
theorem yield
\[
\mathcal E^f_{S,\tau}(Y_\tau)
\leq
\mathcal E^f_{S,\tau}(D_\tau).
\]
Since $D$ is a strong $\mathcal E^f$-martingale,
\[
\mathcal E^f_{S,\tau}(D_\tau)=D_S.
\]
Consequently,
\[
R_S
=
\operatorname*{ess\,sup}_{\tau\in\mathcal T_{S,T}}
\mathcal E^f_{S,\tau}(Y_\tau)
\leq D_S.
\]
Taking the essential infimum over
$D\in\mathfrak C_S^f(Y)$, we obtain
\begin{equation}
\label{eq:dual-lower-bound}
R_S
\leq
\operatorname*{ess\,inf}_{D\in\mathfrak C_S^f(Y)}
D_S.
\end{equation}

We now construct a strong $\mathcal E^f$-martingale
which dominates $Y$ and attains the above bound.

Fix $S\in\mathcal T_{0,T}$ and consider
the following initial-value stochastic equation
on $[S,T]$:
\begin{equation}
\label{eq:optimal-domineering-S}
\left\{
\begin{aligned}
-dD_t^{*,S}
={}&
f(t,D_t^{*,S},Z_t,k_t)\,dt
-Z_t\,dW_t\\
&-
\int_E k_t(e)\,\widetilde N(dt,de),\\
D_S^{*,S}
={}&R_S.
\end{aligned}
\right.
\end{equation}

Since $Z$ and $k$ are fixed processes belonging,
respectively, to $\mathbb H^2$ and
$\mathbb H_\nu^2$, and since $f$ is uniformly
Lipschitz continuous with respect to its
$y$-variable, the above initial-value
stochastic equation admits a unique adapted
solution $D^{*,S}\in\mathcal S^2$ on $[S,T]$.

For any stopping times $\sigma,\tau$ satisfying
$S\leq\sigma\leq\tau\leq T$, we have
\[
\begin{aligned}
D_\sigma^{*,S}
={}&
D_\tau^{*,S}
+\int_\sigma^\tau
f(t,D_t^{*,S},Z_t,k_t)\,dt\\
&-\int_\sigma^\tau Z_t\,dW_t
-\int_\sigma^\tau\int_E
k_t(e)\,\widetilde N(dt,de).
\end{aligned}
\]
By uniqueness of the solution of the BSDE,
\[
D_\sigma^{*,S}
=
\mathcal E^f_{\sigma,\tau}(D_\tau^{*,S}).
\]
Thus, $D^{*,S}$ is a strong
$\mathcal E^f$-martingale on $[S,T]$.

It remains to prove that $D^{*,S}$ dominates
the payoff process $Y$.

Set
\[
\Delta_t^S:=D_t^{*,S}-R_t,
\qquad S\leq t\leq T.
\]
Subtracting the reflected BSDE for $R$
from equation~\eqref{eq:optimal-domineering-S},
we obtain
\[
d\Delta_t^S
=
-\bigl(
f(t,D_t^{*,S},Z_t,k_t)
-
f(t,R_t,Z_t,k_t)
\bigr)\,dt
+dA_t+dB_{t-}.
\]

Define
\[
a_t^S
:=
\frac{
f(t,D_t^{*,S},Z_t,k_t)
-
f(t,R_t,Z_t,k_t)
}{
D_t^{*,S}-R_t
}
\mathbf 1_{\{D_t^{*,S}\neq R_t\}}.
\]
By the Lipschitz continuity of $f$,
\[
|a_t^S|\leq K.
\]
Hence,
\begin{equation}
\label{eq:delta-S}
d\Delta_t^S
=
-a_t^S\Delta_t^S\,dt
+dA_t+dB_{t-},
\qquad
\Delta_S^S=0.
\end{equation}

Consider the integrating factor
\[
\Gamma_t^S
:=
\exp\left(
\int_S^t a_u^S\,du
\right),
\qquad S\leq t\leq T.
\]
Then $\Gamma^S$ is continuous, strictly
positive and of finite variation, with
\[
d\Gamma_t^S
=
a_t^S\Gamma_t^S\,dt,
\qquad
\Gamma_S^S=1.
\]

Applying the Gal'chouk--Lenglart formula
(see Gal'chuk~\cite{Galchuk1981})
to the product $\Gamma^S\Delta^S$,
and using the continuity and finite-variation
property of $\Gamma^S$, we obtain
\[
\begin{aligned}
d(\Gamma_t^S\Delta_t^S)
={}&
\Gamma_t^S\,d\Delta_t^S
+
\Delta_t^S\,d\Gamma_t^S\\
={}&
\Gamma_t^S\,dA_t
+
\Gamma_t^S\,dC_{t-}.
\end{aligned}
\]
Indeed, the drift terms cancel because
\[
-\Gamma_t^S a_t^S\Delta_t^S\,dt
+
a_t^S\Gamma_t^S\Delta_t^S\,dt
=0.
\]

Integrating over the stochastic interval
$[S,t]$ and using $\Delta_S^S=0$, we obtain
\begin{equation}
\label{eq:integrating-factor-S}
\Gamma_t^S\Delta_t^S
=
\int_{(S,t]}\Gamma_u^S\,dA_u
+
\int_{[S,t)}\Gamma_u^S\,dC_u,
\qquad S\leq t\leq T.
\end{equation}

Here, the first integral corresponds to the
right-continuous increasing process $A$,
whereas the second integral corresponds
to the left-continuous process $B_{-}$.
In particular, the interval $[S,t)$ includes
the possible jump of $C$ at $S$, which acts
through the right jump of $R$ immediately
after $S$, while excluding the jump at $t$.

Since $A$ and $C$ are nondecreasing and
$\Gamma^S$ is strictly positive,
equation~\eqref{eq:integrating-factor-S} yields
\[
\begin{aligned}
D_t^{*,S}-R_t
={}&
(\Gamma_t^S)^{-1}
\left(
\int_{(S,t]}\Gamma_u^S\,dA_u
+
\int_{[S,t)}\Gamma_u^S\,dC_u
\right)\\
\geq{}&0,
\qquad S\leq t\leq T.
\end{aligned}
\]

Therefore,
\[
D_t^{*,S}\geq R_t\geq Y_t,
\qquad S\leq t\leq T,
\]
up to the corresponding almost-sure
interpretation at stopping times.

Consequently,
\[
D^{*,S}\in\mathfrak C_S^f(Y).
\]

Since $D_S^{*,S}=R_S$, it follows that
\begin{equation}
\label{eq:dual-upper-bound-S}
\operatorname*{ess\,inf}_{D\in\mathfrak C_S^f(Y)}
D_S
\leq
D_S^{*,S}
=
R_S.
\end{equation}

Combining
\eqref{eq:dual-lower-bound}
and
\eqref{eq:dual-upper-bound-S},
we conclude that
\[
\boxed{
\operatorname*{ess\,sup}_{\tau\in\mathcal T_{S,T}}
\mathcal E^f_{S,\tau}(Y_\tau)
=
R_S
=
\operatorname*{ess\,inf}_{D\in\mathfrak C_S^f(Y)}
D_S.
}
\]

Moreover, the essential infimum is attained
by $D^{*,S}$. Hence $D^{*,S}$ is an optimal
$\mathcal E^f$-domineering claim on $[S,T]$.
\end{proof}

The following proposition characterizes the attainment of the global dual representation in terms of the martingale property of the nonlinear Snell envelope and the absence of reflection. In particular, it provides necessary and sufficient conditions under which a global domineering claim coincides with the nonlinear Snell envelope at a given stopping time.
\begin{proposition}
\label{prop:global-dual-attainment}

Let $Y\in\mathcal S^2$ be a left-limited and
right-upper semicontinuous optional process,
and let $f$ be a Lipschitz driver satisfying
Assumption~A.

Let $(R,Z,k,A,B)$ be the solution of the
reflected BSDE associated with $(f,Y)$.

Denote by $\mathfrak C^f(Y)$ the global class
of all strong $\mathcal E^f$-martingales
$D\in\mathcal S^2$ on $[0,T]$ satisfying
\[
D_\tau\geq Y_\tau
\quad\mathbb P\text{-a.s. for every }
\tau\in\mathcal T_{0,T}.
\]

Then, for every $S\in\mathcal T_{0,T}$,
the following assertions are equivalent:

\begin{enumerate}[label=\textnormal{(\roman*)}]

\item There exists $D\in\mathfrak C^f(Y)$
such that
\[
D_S=R_S
\quad\mathbb P\text{-a.s.}
\]

\item The nonlinear Snell envelope $R$
is a strong $\mathcal E^f$-martingale
on the stochastic interval $[0,S]$.

\item The reflection processes satisfy
\[
A_S=B_{S-}=0
\quad\mathbb P\text{-a.s.}
\]

\end{enumerate}

Consequently, the global dual representation
attains the nonlinear Snell envelope at $S$
if and only if no reflection has occurred
before or at $S$, in the sense of (iii).

\end{proposition}

\begin{proof}

We prove the implications
$(i)\Rightarrow(ii)\Rightarrow(iii)
\Rightarrow(i)$.

\medskip
\noindent
\textbf{Step 1: $(i)\Rightarrow(ii)$.}

Suppose that there exists
$D\in\mathfrak C^f(Y)$ such that
\[
D_S=R_S
\quad\mathbb P\text{-a.s.}
\]

Since $D$ dominates $Y$ and is a strong
$\mathcal E^f$-martingale, the
reflected-BSDE characterization yields
\[
R_\sigma\leq D_\sigma
\]
for every $\sigma\in\mathcal T_{0,T}$.

Let $\sigma\in\mathcal T_{0,T}$ satisfy
$\sigma\leq S$.

By the strong $\mathcal E^f$-martingale
property of $D$, we have
\[
D_\sigma
=
\mathcal E^f_{\sigma,S}(D_S)
=
\mathcal E^f_{\sigma,S}(R_S).
\]

On the other hand, since $R$ is a strong
$\mathcal E^f$-supermartingale,
\[
\mathcal E^f_{\sigma,S}(R_S)
\leq R_\sigma.
\]

Combining these inequalities gives
\[
R_\sigma
\leq D_\sigma
=
\mathcal E^f_{\sigma,S}(R_S)
\leq R_\sigma.
\]

Therefore,
\[
R_\sigma
=
\mathcal E^f_{\sigma,S}(R_S)
\quad\mathbb P\text{-a.s.}
\]

By time consistency, for any stopping
times $\sigma\leq\tau\leq S$,
\[
R_\sigma
=
\mathcal E^f_{\sigma,\tau}(R_\tau).
\]

Thus $R$ is a strong
$\mathcal E^f$-martingale on $[0,S]$.

\medskip
\noindent
\textbf{Step 2: $(ii)\Rightarrow(iii)$.}

Suppose that $R$ is a strong
$\mathcal E^f$-martingale on $[0,S]$.

By the reflected BSDE, for
$t\leq S$,
\[
\begin{aligned}
R_t={}&R_S
+\int_t^S f(u,R_u,Z_u,k_u)\,du\\
&-\int_t^S Z_u\,dW_u\\
&-\int_t^S\int_E
k_u(e)\,\widetilde N(du,de)\\
&+A_S-A_t+B_{S-}-B_{t-}.
\end{aligned}
\]

Since $R$ is also a strong
$\mathcal E^f$-martingale on $[0,S]$,
it admits a BSDE representation without
an additional finite-variation term.

By uniqueness of the canonical
semimartingale decomposition on the
stochastic interval $[0,S]$,
the increasing reflection term must vanish:
\[
A_t+B_{t-}=0,
\qquad 0\leq t\leq S.
\]

In particular,
\[
A_S+B_{S-}=0
\quad\mathbb P\text{-a.s.}
\]

Since $A$ and $B$ are nondecreasing,
with $A_0=B_{0-}=0$, we conclude that
\[
A_S=B_{S-}=0
\quad\mathbb P\text{-a.s.}
\]

\medskip
\noindent
\textbf{Step 3: $(iii)\Rightarrow(i)$.}

Suppose that
\[
A_S=B_{S-}=0
\quad\mathbb P\text{-a.s.}
\]

Consider the initial-value stochastic
equation on $[0,T]$:
\[
\left\{
\begin{aligned}
-dD_t^*
={}&f(t,D_t^*,Z_t,k_t)\,dt
-Z_t\,dW_t\\
&-\int_E k_t(e)\,
\widetilde N(dt,de),\\
D_0^*={}&R_0.
\end{aligned}
\right.
\]

As in the proof of
Proposition~\ref{prop:dual-representation},
this equation admits a unique solution
$D^*\in\mathcal S^2$.

Moreover, $D^*$ is a strong
$\mathcal E^f$-martingale on $[0,T]$.

Set
\[
\Delta_t:=D_t^*-R_t.
\]

Subtracting the reflected BSDE from
the equation for $D^*$ yields
\[
d\Delta_t
=
-a_t\Delta_t\,dt
+dA_t+dB_{t-},
\qquad \Delta_0=0,
\]
where
\[
a_t=
\frac{
f(t,D_t^*,Z_t,k_t)
-f(t,R_t,Z_t,k_t)
}{
D_t^*-R_t
}
\mathbf 1_{\{D_t^*\neq R_t\}}.
\]

The Lipschitz property implies
$|a_t|\leq K$.

Define
\[
\Gamma_t
=
\exp\left(\int_0^t a_u\,du\right).
\]

Applying the Gal'chouk--Lenglart
integration-by-parts formula gives
\[
\Gamma_t\Delta_t
=
\int_{(0,t]}\Gamma_u\,dA_u
+
\int_{[0,t)}\Gamma_u\,dB_u.
\]

Since $\Gamma$ is strictly positive
and $A,B$ are nondecreasing,
\[
D_t^*\geq R_t\geq Y_t
\]
at every stopping time.

Consequently,
\[
D^*\in\mathfrak C^f(Y).
\]

Furthermore, using
$A_S=B_{S-}=0$, we obtain
\[
\begin{aligned}
\Gamma_S(D_S^*-R_S)
&=
\int_{(0,S]}\Gamma_u\,dA_u\\
&\quad+
\int_{[0,S)}\Gamma_u\,dB_u\\
&=0.
\end{aligned}
\]

Since $\Gamma_S>0$,
\[
D_S^*=R_S
\quad\mathbb P\text{-a.s.}
\]

Thus assertion (i) holds.

The proof is complete.
\end{proof}

We now establish a necessary contact condition
for an optimal domineering claim. More precisely,
under the strict comparison property, we show
that any optimal stopping time must be a contact
time between the payoff process and a global
domineering claim attaining the optimal value.

\begin{proposition}[Necessary contact condition for optimality]
\label{prop:optimal-contact-claim}

Let $Y\in\mathcal S^2$ be an optional payoff
process, and let $C\in\mathfrak C^f(Y)$
be a global $\mathcal E^f$-domineering claim.

Assume that $C$ attains the nonlinear
optimal stopping value at time zero, i.e.,
\[
C_0=R_0
=
\sup_{\tau\in\mathcal T_{0,T}}
\mathcal E^f_{0,\tau}(Y_\tau).
\]

Suppose that the optimal stopping problem
admits an optimal stopping time
$\tau^*\in\mathcal T_{0,T}$.

If $\mathcal E^f_{0,\tau^*}$ satisfies
the strict comparison property, then
\[
C_{\tau^*}=Y_{\tau^*}
\quad\mathbb P\text{-a.s.}
\]

Consequently,
\[
\operatorname*{ess\,sup}_{\tau\in\mathcal T_{0,T}}
(Y_\tau-C_\tau)=0
\quad\mathbb P\text{-a.s.}
\]

\end{proposition}

\begin{proof}

Since $\tau^*$ is optimal and
$C_0=R_0$, we have
\[
C_0
=
R_0
=
\mathcal E^f_{0,\tau^*}
(Y_{\tau^*}).
\]

Moreover, since $C$ is a strong
$\mathcal E^f$-martingale,
\[
C_0
=
\mathcal E^f_{0,\tau^*}
(C_{\tau^*}).
\]

Therefore,
\[
\mathcal E^f_{0,\tau^*}
(Y_{\tau^*})
=
\mathcal E^f_{0,\tau^*}
(C_{\tau^*}).
\]

Since $C$ dominates $Y$,
\[
Y_{\tau^*}\leq C_{\tau^*}
\quad\mathbb P\text{-a.s.}
\]

By the strict comparison property,
\[
Y_{\tau^*}=C_{\tau^*}
\quad\mathbb P\text{-a.s.}
\]

On the other hand, the domination
property gives
\[
Y_\tau-C_\tau\leq 0
\quad\mathbb P\text{-a.s.}
\]
for every $\tau\in\mathcal T_{0,T}$.

Hence,
\[
\begin{aligned}
0
&=Y_{\tau^*}-C_{\tau^*}\\
&\leq
\operatorname*{ess\,sup}_{\tau\in\mathcal T_{0,T}}
(Y_\tau-C_\tau)\\
&\leq 0.
\end{aligned}
\]

Consequently,
\[
\boxed{
\operatorname*{ess\,sup}_{\tau\in\mathcal T_{0,T}}
(Y_\tau-C_\tau)=0
}
\quad\mathbb P\text{-a.s.}
\]

This completes the proof.
\end{proof}

\section*{Acknowledgments}

Oussama LAMHADOUL acknowledges the financial support
of the National Center for Scientific and Technical
Research (CNRST), Morocco, through the
PhD-ASsociate Scholarship -- PASS Program.

\end{document}